# Large-scale rail freight flow assignment with the tree-shaped path

Bo-Liang Lin[*]

*School of Traffic and Transportation, Beijing Jiaotong University, Beijing 100044, China*

**Abstract**: The rail freight flow assignment (RFFA) is a key point of network design and operations management. Unlike a road network, the freight flow assignment pattern in a rail system has unique characteristics, and the most significant one is the tree-shaped path. More specifically, when two or more shipments arrive at a rail yard, those intended for the same destination will be treated as a single freight flow and be shipped through the same path during the remainder of the trip, regardless of their origins. This paper reveals the forming mechanism of the tree-shaped path from the perspective of rail yard operations, and develops a nonlinear programming model for the RFFA problem. Considering the notably large number of quadratic terms generated from the tree-shaped constraints, the model is difficult to solve even for small-size instances. We thus introduce binary decision variables and modify the model into a nonlinear binary programming problem. In order to solve large-scale RFFA problems, we outline a simulated annealing algorithm.

## 1. Introduction

It is known that the congestion phenomena that frequently occurs in the urban road system rarely happens in the rail network, primarily because the railway system operates in accordance with train schedules. The basis of train schedule design is the rail freight or passenger flow assignment, especially the rail freight flow assignment (RFFA). The RFFA is not only the foundation of train formation planning and train scheduling, but also the basis of rail network design.

The key to solving the RFFA problem is the optimization of the rail freight car path. Unlike the road network, the rail freight car path has its own unique characteristics. The train drivers have no right to choose or change the train route freely and the change of running direction of a train is effected by moving the position of the railroad switches, which used to be completed by ground switchmen. However, the switchmen actually do nothing but executing commands given by the dispatching center. Tracing further up the command chain, we can find that the commands from the dispatching center are actually specified in a train path plan, which is set in advance. Because of the capacity shortage of some links in a rail network, only a portion of origin-destination (O-D) pairs can be shipped through their shortest paths. The path plan, aiming at minimizing the total transport cost in the rail network, exactly determines which O-D pairs should be assigned to their shortest paths and which O-D pairs should be selected to make a detour. This is a typical network flow assignment problem with the system optimal principle.

Taking the China railway system as an example, according to the data presented by the National Railway Administration (NRA) of China, by 2017, China had built a railway network of 127 thousand kilometers. Approximately 3,000 stations offer transport services for hundreds of thousands of freight O-D pairs (shipments) in this network. The China railway has become one of the largest and busiest railway systems in the world, and therefore, solving the RFFA problem for such a large-scale network is theoretically challenging and practically tricky.

---
[*] E-mail address: bllin@bjtu.edu.cn (B.-L. Lin).

The RFFA is usually involved in the railway network design process (see, e.g., Lin et al., 2017). Kuby et al. (2001) presented a rail network design model in a project supported by the World Bank and China's Ministry of Railways. The model removed all small traffic demands below a certain cutoff level from the O-D matrix and preloaded them onto their shortest paths, and this approach yielded great benefits for model size and performance. However, in China railway system, small traffic demands usually have longer travel distances than larger ones. The consequence of adopting this model means that long-distance freight flows go along their shortest paths while those with short distances might have to detour. This result did not mimic real railway operations in China, which prefer to allow long-distance shipments to take detours. Lin et al. (2002) addressed the optimization decision for railroad network design problems to determine where a new railroad line should be constructed or an old railroad should be upgraded from the viewpoint of project investment and cost of flow assignment. Tian et al. (2011) developed an arc-path model for the car flow assignment problem and solved infeasible flows by setting virtual arcs. In addition, the sets of available paths were created by limiting the detour ratio of the paths.

## 2. Problem description for the tree-shaped path

Unlike the road traffic assignment, in most railway operations management, when two or more freight flows to the same destination merge at a railway yard, they will be considered as one flow and run on the same path during the remainder of the trip, even they may come from different origins. This consideration forms a unique pattern of the flow paths at the macro level, which we define as the tree-shaped path.

The real-life background of the tree-shaped path can be nicely illustrated by a simple example. A portion of the railway network in the north and west of China is shown in Fig. 1, and selected train service plans are listed as follows.

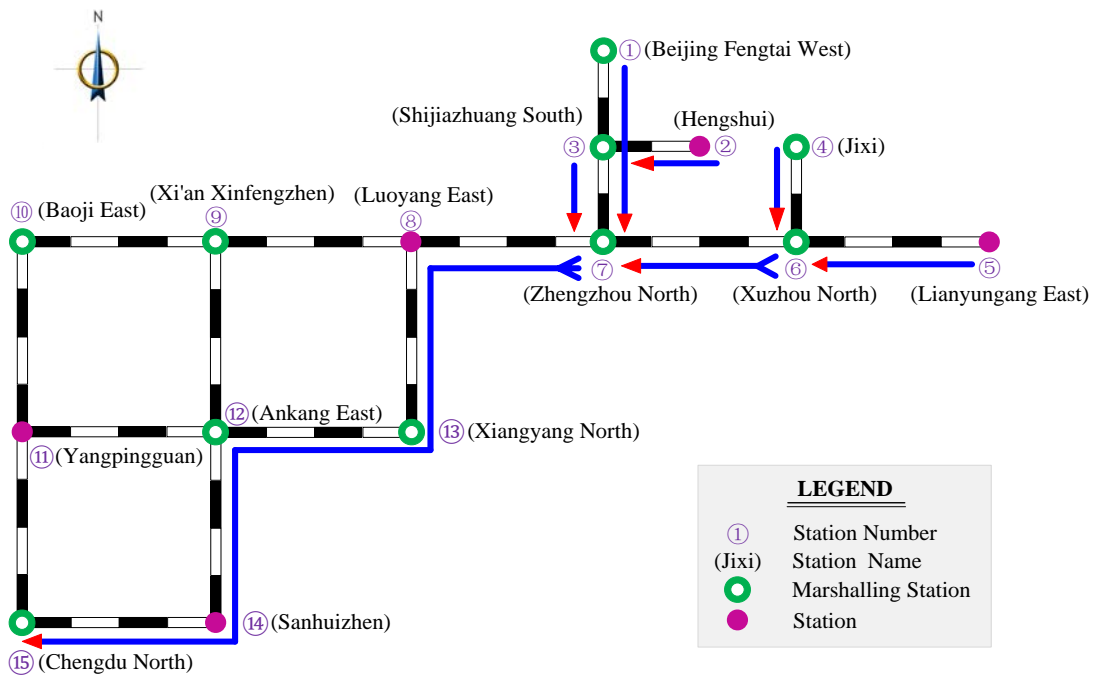

**Fig. 1.** Portions of the China railway network.

The cars of direct freight train ①→⑦ are most come from North China, and their destinations are distributed in the south and west of China and the region of ⑦ (Zhengzhou North marshalling station, see Fig. 2). When the train arrives at ⑦, it will be received into the yard XD. After the arrival operation is completed, the train will be pushed over the hump. The cars destined for the southwest regions of China, such as ⑮ and beyond, are rolled onto a certain track (e.g., Track 7) in the yard XB. The cars heading to western regions are rolled onto another track (e.g., Track 2). The cars destined for southern regions and local cars are rolled onto two other tracks (e.g., Track 4 and Track 1). These operations are similar to those of the train ③→⑦. When this train arrives at ⑦, the cars are rolled onto the corresponding tracks according to their destinations. For example, cars heading to the southwest are rolled to Track 7 in the yard XB, and cars destined for the west, south and the region of ⑦ are rolled to Tracks 2, 4 and 1 in the yard XB respectively. As the trains continually arrive and break up, the number of cars accumulated on each track in the yard XB increases. When the accumulated cars on the track reach the full tonnage, they will be hauled to the yard XF to form a new train. After departure operations, the train will leave ⑦ for its destination along a planned path. For example, the new train ⑦→⑮ (indicated by the red dashed lines in Fig. 2) carrying cars destined for the southwest regions leaves for ⑮ in accordance with the path ⑦→⑬→⑫→⑭→⑮ (see Fig. 1).

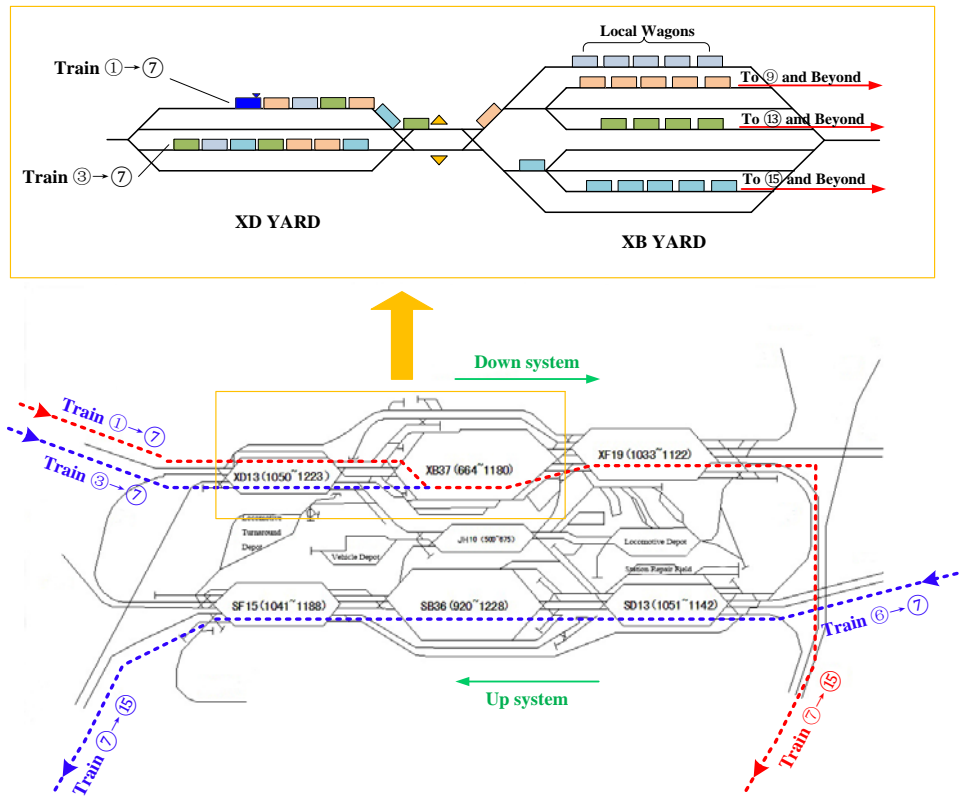

**Fig. 2.** The Zhengzhou North marshalling station.

We can also analyze train ⑥→⑦ (shown in the blue dashed line in Fig. 2). The cars of this train mostly come from East China, such as ④ and ⑤. After this train arrives at ⑦ and is broken up, cars are sorted, and those heading to the southwest are grouped into the train ⑦→⑮ (shown in the blue dashed line in Fig. 2) which bounds for ⑮ in the up system. The train ⑦→⑮ selects the same path of the train ⑦→⑮ mentioned above, i.e., ⑦→⑬→⑫→⑭→⑮, regardless of whether the trains are formed in the up system or in the down system.

In fact, with further analysis, we find that the cars in the train ①→⑦ come from ① and the areas on its north. Similar cases include the trains ③→⑦ and ⑥→⑦. It is observed that the cars carried by the trains ①→⑦, ③→⑦ and ⑥→⑦ are grouped into the train ⑦→⑮ at ⑦ if the destination of these cars is ⑮ and beyond. From a broader perspective, the freight flows from different origins heading to the same destination finally converge into a large flow after several merges. Their path trace resembles a tree (see Fig. 3) with the tree root located at ⑮.

Therefore, we define such a running trace of the rail freight flow as the tree-shaped path.

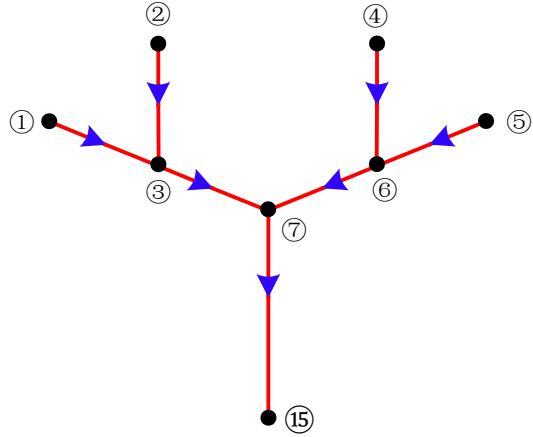

**Fig. 3.** The tree-shaped path.

## 3. Mathematical formulation of the RFFA problem

In this section, we present rail freight flow assignment models based on the tree-shaped path. The models aim at determining: (i) which path is chosen by each O-D pair, (ii) whether the flow volume exceeds its capacity on each link, (iii) whether the physical path of each O-D pair satisfies the tree-shaped path constraint, and (iv) whether the flow conservation equation holds at each node. The objective of the models is to minimize the total transport cost in a rail network.

### 3.1 The nonlinear programming model

The notations used in the mathematical formulations are defined in Table 1.

**Table 1** Notations used in the formulation.

| Symbol | Definition |
|---|---|
| $V$ | Set of all stations in a rail network |
| $V(i)$ | Set of stations adjacent to $i$ |
| $E$ | Set of all links in a rail network |
| $N$ | Set of all O-D pairs |
| $R_{ij}$ | Available capacity of arc $i \rightarrow j$ |
| $c_{ij}$ | Unit transport cost on arc $i \rightarrow j$ |
| $N_{ij}$ | Flow volume of the O-D pair that originates at station $i$ and is destined for station $j$ |
| $f_{ij}$ | The actual flow volume from station $i$ to $j$ (which consists of the original O-D volume $N_{ij}$ and transfer flow volume at station $i$) |
| $f_{ij}^k$ | Decision variables, the volume of selected flow portion of $f_{ij}$ that chooses adjacent station $k$ of $i$ as the first front station in its travel itinerary |

In order to illustrate the mathematical model based on the tree-shaped path, a simplified directed rail network is constructed in Fig. 4.

A general formula which can describe the relationship between $f_{ij}$ and $f_{ij}^k$ from the above two perspectives will be expressed by

$$f_{ij} = \sum_{k \in V(i)} f_{ij}^k \quad \forall i \in V, j \in V, i \neq j \tag{1}$$

which also can be expressed as:

$$f_{ij} = N_{ij} + \sum_{h \in V(i)} f_{hj}^i \quad \forall i \in V, j \in V, i \neq j \tag{2}$$

Conbining Eq. (1) andof Eq. (2), we can write the standard flow conservation equation as follows:

$$\sum_{k \in V(i)} f_{ij}^k - \sum_{h \in V(i)} f_{hj}^i = N_{ij} \quad \forall i \in V, j \in V, i \neq j \tag{3}$$

Under the tree-shaped path constraints, flow $f_{ij}$ can only choose one path. Equivalently, the flow can only choose one adjacent station even if there are several stations available. We use the following formula to describe such a condition:

$$\sum_{k \in V(i)} \sum_{h \in V(i), h \neq k} f_{ij}^k f_{ij}^h = 0 \quad \forall i \in V, j \in V, i \neq j \tag{4}$$

**Remark.** Here we present a proof that Constraint (4) can guarantee the tree-shaped path in a general network using the graph theory. For a non-zero-volume flow $f_{ij}$ ($i, j \in V, i \neq j$), we consider the following two cases: (I) $f_{ij}$ is diverted into serval flows at its first front station. Without loss of generality, we assume two of the flows choose $k$ and $h$ as the first front station, respectively, i.e. $f_{ij}^k > 0$ and $f_{ij}^h > 0$. Thus, we have $\sum_{k \in V(i), h \in V(i), h \neq k} f_{ij}^k f_{ij}^h \geq f_{ij}^k f_{ij}^h > 0$. Clearly, this is contradicted with Constraint (4). (II) $f_{ij}$ is diverted into serval flows at a farer station $s$. Without loss of generality, we assume two of the flows choose $u$ and $v$ as the first front station, respectively, i.e. $f_{sj}^u > 0$ and $f_{sj}^v > 0$. This will result in $\sum_{u \in V(s), v \in V(s), u \neq v} f_{sj}^u f_{sj}^v \geq f_{sj}^u f_{sj}^v > 0$, which is also contradicted with Constraint (4). Therefore, under Constraint (4), both above two cases do not hold and $f_{ij}$ is not allowed to be split. This means that $f_{ij}$ can only choose one path from $i$ to $j$. Given the proposition that "graph $T$ is a tree" is equivalent to "any two vertices of Graph $T$ are connected by exactly one path" (see Wilson, 1996), Constraint (4) can guarantee the tree-shaped path.

To minimize the ton-km or generalized ton-km, we can formulate the RFFA problem as a nonlinear programming model:

(RFFA I) $$\min \sum_{i \in V} \sum_{k \in V(i)} c_{ik} \sum_{j \in V} f_{ij}^k \tag{5}$$

s.t. $$\sum_{k \in V(i)} f_{ij}^k - \sum_{h \in V(i)} f_{hj}^i = N_{ij} \quad \forall i \in V, j \in V, i \neq j \tag{6}$$

$$\sum_{k \in V(i), h \in V(i), h \neq k} f_{ij}^k f_{ij}^h = 0 \quad \forall i \in V, j \in V, i \neq j \tag{7}$$

$$\sum_{j \in V} f_{ij}^k \leq R_{ik} \quad \forall i \in V, k \in V(i) \tag{8}$$

$$f_{ij}^k \geq 0 \quad \forall i \in V, j \in V, k \in V(i), i \neq j \tag{9}$$

In model RFFA I, the objective function is to minimize the total transportation cost. Constraint (6) is the flow conservation equation. Constraint (7) ensures that all flows finish their journeys along a tree-shaped path. Constraint (8) guarantees that the flow volumes on each link are less than the available capacity.

After solving the model RFFA I, we can obtain the path of $i \to j$ by tracing the variable $f_{ij}^k$, whose value is non-zero. For example, if non-zero variables $f_{27}^3, f_{37}^8, f_{87}^6, f_{67}^4, f_{47}^7$ appeared in an optimized solution, the path of the O-D pair $2 \to 7$ should be $2 \to 3 \to 8 \to 6 \to 4 \to 7$.

### 3.2 Infeasible flows

It should be noted that the model RFFA I cannot guarantee obtaining a feasible solution because if certain links in the network are in a fully loaded state, some flows cannot reach their destinations through the links. In other words, it is impossible to assign all O-D pairs because one or several bottlenecks exist in the network. To obtain feasible solutions, some O-D pairs must be abandoned from the perspective of system optimization. There is no doubt that we wish the potential transportation loss of the O-D pairs is minimal. Let $U_{ij}$ represent the portion of demand $N_{ij}$ that cannot be met, and let $P_{ij}$ denote the loss of income (shadow price) when a unit of transportation demand from $i$ to $j$ is given up. In this way, we can modify RFFA I as follows:

(RFFA II)
$$\min \sum_{i \in V} \sum_{k \in V(i)} c_{ik} \sum_{j \in V} f_{ij}^k + \sum_{(i,j) \in N} P_{ij} U_{ij} \tag{10}$$

s.t.
$$\sum_{k \in V(i)} f_{ij}^k - \sum_{h \in V(i)} f_{hj}^i + U_{ij} = N_{ij} \quad \forall i \in V, j \in V, i \neq j \tag{11}$$

$$\sum_{k \in V(i), h \in V(i), h \neq k} f_{ij}^k f_{ij}^h = 0 \quad \forall i \in V, j \in V, i \neq j \tag{12}$$

$$\sum_{j \in V} f_{ij}^k \leq R_{ik} \quad \forall i \in V, k \in V(i) \tag{13}$$

$$0 \leq U_{ij} \leq N_{ij} \quad \forall \{(i,j) : N_{ij} \in N\} \tag{14}$$

$$f_{ij}^k \geq 0 \quad \forall i \in V, j \in V, k \in V(i), i \neq j \tag{15}$$

The difference between RFFA II and RFFA I is embodied in the objective function and the flow conservation equation. In model RFFA II, on the one hand, we add the losses of transportation revenue to the objective function that aims to minimize $U_{ij}$. On the other hand, we also add $U_{ij}$ to the flow conservation equation. Consequently, as $U_{ij}$ decreases, the feasible flow volumes through the relevant links will increase, however, it will be ultimately limited by the link capacity Constraint (13). Thus, under the combined effect of the objective function and constraint conditions, we can obtain a reasonable set of infeasible flows. In addition, Constraint (14) is added to ensure that the infeasible flow is less than the original O-D demand. In other words, the worst scenario is that the volume of a certain O-D pair is given up entirely rather than partly.

In fact, the value of $U_{ij}$ can also be determined by virtual arcs, i.e., for a given network $(V, E)$ and a set of O-D pairs $N$, we can add a virtual arc $e_{ij}$ with a general mileage of $P_{ij}$ for all node pairs $\{(i, j) : N_{ij} \in N\}$ into the network. If there is already an arc between $i$ and $j$ in the original network, we can add a virtual node between them, e.g., node $k$ (of course, it must be different from both the existing nodes and virtual nodes). In this way, the virtual corridor between $i$ and $j$ consists of virtual arcs $e_{ik}$ and $e_{kj}$ (with an equal cost of $P_{ij}/2$). Naturally, the corresponding virtual path from $i$ to $j$ is $i \rightarrow k \rightarrow j$. Let $\Delta E$ and $\Delta V$ respectively denote the sets of all virtual arcs and virtual nodes, then in the extended network $(V + \Delta V, E + \Delta E)$ (see Fig. 4), model RFFA II will regress to model RFFA I. It is worth mentioning that the value of $P_{ij}$ should be larger than the length of the possible longest path. A reasonable method for valuing $P_{ij}$ is to adopt the sum of all arcs' weights.

### 3.3 The nonlinear binary programming model

The above two models contain many quadratic terms due to the existence of tree-shaped path constraints. In fact, Eq. (4) can be equally expressed by:

$$\sum_{k \in V(i)} I(f_{ij}^k) = 1 \quad \forall i \in V, j \in V, i \neq j \tag{16}$$

where the indicator function $I(x)$ is defined as follows:

$$I(x) = \begin{cases} 1 & \text{If } x > 0 \\ 0 & \text{Otherwise} \end{cases} \tag{17}$$

For Fig. 4, Eq. (16) can be expanded as

$$I(f_{ij}^2) + I(f_{ij}^3) + I(f_{ij}^4) = 1 \tag{18}$$

In fact, the indicator function $I(x)$ can also be replaced by 0-1 decision variables. We define the freight flow path choice decision variables as follows:

$$x_{ij}^k = \begin{cases} 1 & \text{If } f_{ij} \text{ choose the adjacent station } k \text{ as the first front station} \\ 0 & \text{Otherwise} \end{cases} \quad (19)$$

Therefore, the flow conservation Eq. (2) can be expressed by

$$f_{ij} - \sum_{s \in V(i)} f_{sj} x_{sj}^i = N_{ij} \quad \forall i \in V, j \in V, i \neq j \quad (20)$$

In this way, the model RFFA I can be reformulated in the following form:

(RFFA III) 
$$\min \sum_{i \in V} \sum_{k \in V(i)} c_{ik} \sum_{j \in V} f_{ij} x_{ij}^k \quad (21)$$

s.t. 
$$\sum_{k \in V(i)} x_{ij}^k = 1 \quad \forall i \in V, j \in V, i \neq j \quad (22)$$

$$\sum_{j \in V} f_{ij} x_{ij}^k \leq R_{ik} \quad \forall i \in V, k \in V(i) \quad (23)$$

$$x_{ij}^k \in \{0,1\} \quad \forall i \in V, j \in V, k \in V(i), i \neq j \quad (24)$$

where flow $f_{ij}$ can be derived by the following recursive expression:

$$f_{ij} = N_{ij} + \sum_{s \in V(i)} f_{sj} x_{sj}^i \quad \forall i \in V, j \in V, i \neq j \quad (25)$$

Of course, the model RFFA II can also be reformulated using the 0-1 decision variables $x_{ij}^k$ similarly:

(RFFA IV) 
$$\min \sum_{i \in V} \sum_{k \in V(i)} c_{ik} \sum_{j \in V} f_{ij} x_{ij}^k + \sum_{(i,j) \in N} P_{ij} U_{ij} \quad (26)$$

s.t. 
$$f_{ij} - \sum_{s \in V(i)} f_{sj} x_{sj}^i + U_{ij} = N_{ij} \quad \forall i \in V, j \in V, i \neq j \quad (27)$$

$$\sum_{k \in V(i)} x_{ij}^k = 1 \quad \forall i \in V, j \in V, i \neq j \quad (28)$$

$$\sum_{j \in V} f_{ij} x_{ij}^k \leq R_{ik} \quad \forall i \in V, k \in V(i) \quad (29)$$

$$0 \leq U_{ij} \leq N_{ij} \quad \forall \{(i,j) : N_{ij} \in N\} \quad (30)$$

$$x_{ij}^k \in \{0,1\} \quad \forall i \in V, j \in V, k \in V(i), i \neq j \quad (31)$$

$$f_{ij} \geq 0 \quad \forall i \in V, j \in V, i \neq j \quad (32)$$

### 3.4 Complexity analysis and the detour ratio

The considered RFFA problem can be viewed as a generalization of the well-known Integer Multicommodity Flow Problem (IMFP). The IMFP has been proved to be NP-complete, for example, in Even et al. (1975). Therefore, the RFFA problem also belongs to the NP-complete class. Note that the RFFA problem is more complicated and computationally more challenging due to its additional tree-shaped path constraints.

For a rail network consisting of $n$ nodes, if the average degree of each node (i.e., number of adjacent links) is $m$, considering that $f_{ij}$ are intermediate variables, the number of decision variables should be $n \times (n-1) \times m$. The number of flow conservation constraints and tree-shaped path constraints are both $n \times (n-1)$ while the number of capacity constraints is $n \times m$. It is worth pointing out that if the degree of node $i$ is $d$, the number of first front stations that can be chosen by $f_{ij}$ is usually less than $d$. In fact, this situation depends on the allowable detour ratio of the freight flow path. When flow $f_{ij}$ chooses $k$ as the first front station, then the minimum path length $\rho_{ij}^k$ will be $\rho_{ij}^k = c_{ik} + \rho_{kj}$, where $\rho_{ij}$ is the shortest path from $i$ to $j$. Therefore, the path detour ratio $\varepsilon$ can be defined as follows:

$$\varepsilon = \rho_{ij}^k / \rho_{ij} \quad (33)$$

Considering that corresponding paths of $\rho_{ij}$ and $\rho_{ij}^k$ might contain common links, e.g., in Fig. 4 the path $i \to 2 \to 8 \to 6 \to 7 \to j$ and path $i \to 4 \to 5 \to 6 \to 7 \to j$ contain two common arcs $6 \to 7$ and $7 \to j$, a more reasonable approach to obtaining $\varepsilon$ is to cut out the common arcs. Let $\Delta$ denote the length of common segments of these two paths, and the relative path detour ratio can be defined as:

$$\varepsilon = (\rho_{ij}^k - \Delta) / (\rho_{ij} - \Delta) \qquad (34)$$

If $\varepsilon$ is greater than a predefined threshold, such as $\varepsilon = 1.4$, we consider that it is unreasonable for $f_{ij}$ to choose station $k$ as the first front station. Obviously, it is unreasonable to choose $k$ as the first front station that makes $\rho_{ij} = \Delta$. For example, if we choose ① as the first front station by flow $f_{ij}$, equation $\rho_{ij} = \Delta$ will be true (see Fig. 4). In this case, the denominator of formula (27) is equal to zero.

## 4. Solving the RFFA using simulated annealing

The SA method was independently described by Kirkpatrick et al. (1983) and Černý (1985). It is an adaptation of the Metropolis-Hastings algorithm, a Monte Carlo method to generate sample states of a thermodynamic system, published by Metropolis et al. (1953).

The motivations of using the SA method are threefold: (1) As early as 1990s, our research group has tried to solve the RFFA problem and the rail freight train connection service network problem using other heuristics (e.g. the genetic algorithm, see Li, 2000). However, it turned to be that the performances of other heuristics were not much better than the SA. (2) The SA algorithm has the advantages of excellent exploration ability (i.e. avoiding local optimums) and the relative ease of implementation, compared to other meta-heuristics, such as the genetic algorithm and the neural network algorithm (see, e.g., Selim and Alsultan, 1991). (3) Inspired by our previous work in successfully solving the problem of large-scale rail freight train connection service network using the SA algorithm (see Lin et al., 2012). To apply the SA method to our problem, we need to modify the model RFFA III first. In the model, the capacity constraint of each arc is a "difficult" constraint, which means that it is not easy to be met for any solution. A widely adopted technique to deal with these constraints is known as the penalty function method, i.e., converting the difficult constraints into the objective function. In this way, we introduce the penalty function for model RFFA III as follows:

$$H(X) = \sum_{i \in V} \sum_{k \in V(i)} \max\{0, \sum_{j \in V} f_{ij} x_{ij}^k - R_{ik}\} \qquad (35)$$

After adding the penalty function into the objective function, the new objective function can be redefined as follows:

$$Z(X) = \sum_{i \in V} \sum_{k \in V(i)} c_{ik} \sum_{j \in V} f_{ij} x_{ij}^k + \lambda H(X) \qquad (36)$$

where $Z(X)$ is referred to as the energy function, and $\lambda$ is a positive penalty parameter. According to our computational experiences, $\lambda$ is set as 600 in this paper.

In theory, SA is able to converge to (near-) optimal solutions for any given finite problem (Granville et al., 1994). In other words, no matter how to select the initial solution, there is little influence on the final solution. Therefore, we can generate an initial solution as follows. For each node pair $i \to j$, we randomly select one variable from the decision variable set $\{x_{ij}^k\}$ and set it to 1 and set all others to 0. In this way, an initial solution is generated and satisfies Constraint (22). Naturally, the flow $f_{ij}$ can be calculated.

For the neighborhood solution, we can generate it as follows. First, we randomly select a node pair $i \to j$ in the corresponding variable set $\{x_{ij}^k\}$. Second, in this set, for the variable that is equal to 1, we set it to 0. Then we randomly select one from the remaining variables in {$x_{ij}^k$} and set it to 1. In this manner, a new neighborhood solution is generated. Finally, this new solution is accepted if it is better than the old solution; otherwise, it will be accepted through the Metropolis rule.

The SA algorithm framework for solving the RFFA problem is shown in Fig. 9. As presented in the figure, the input of the algorithm includes network data (adjacency matrix, link capacity and link length) and shipment data (O-D matrix). In order to reduce the solution space, we can pre-process the decision variables by removing unreasonable choices of adjacent nodes for each shipment based on a given detour ratio $\varepsilon$ (in our algorithm, $\varepsilon = 1.4$). Subsequently, we can perform the general steps of the SA algorithm. The output of the algorithm is the path for each shipment. On basis of the path data, we can calculate the total transport cost, average transport distance for the shipments and the freight volume on each link. Furthermore, we have also developed a software for rail freight assignment that can visualize the output data with a user interface. The software has been a useful tool for precisely analyzing the network flow.

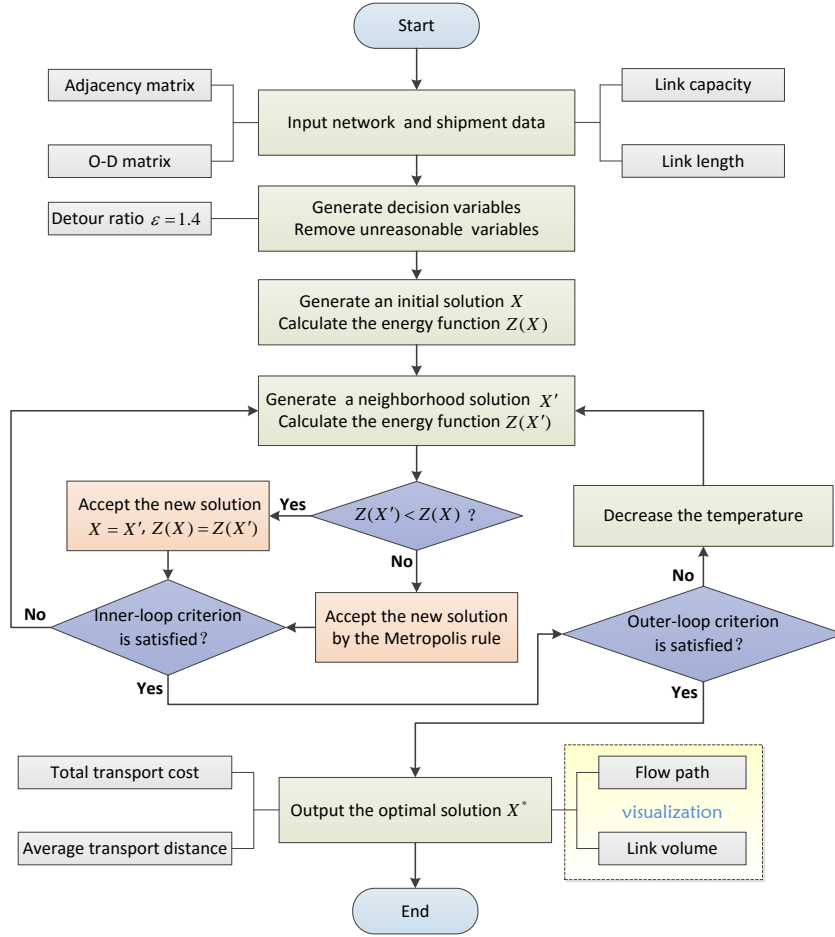

**Fig. 4.** Flow chart of the proposed algorithm.

Under the general framework of the SA algorithm shown in Fig. 4, we make the following major innovations to further improve the performance of the original algorithm. Firstly, instead of setting the total number of inner-loop iterations as a constant (which is a widely used method in the existing literature), we suggest a flexible length $L_T$ of each Markov chain by $L_T = K|\Omega_i|$, where $K$ is set between 3 and 6, and $\Omega_i$ represents the set of neighborhood solutions of the current solution $X_i$. This means that, at each temperature, when the number of generated neighborhood solutions is six times more than $|\Omega_i|$, or when the number of accepted neighborhood solutions is three times more than $|\Omega_i|$, the inner-loop termination criterion will be satisfied. We benefit from this flexile design of Markov chain that at high temperatures, there are not too much solutions to be accepted (because the acceptance rate is relatively high); while at low temperatures, there are not too much solutions to be generated (because the acceptance rate is relatively low). Secondly, to avoid that the energy value decreases too slowly at the early stage of the SA run, we design a piecewise temperature update function (Lin et al., 2012), which can be expressed as follows:

$$T_{m+1} = \begin{cases} T_m \left[ 1 + \dfrac{T_m \ln(1+\delta)^{-1}}{3\sigma(T_m)} \right] & m \leq n \\ \alpha T_m & m > n \end{cases} \quad (37)$$

In the equation, $T_m$ means the temperature at the $m$th iteration, $\delta$, $\alpha$ and $n$ are two pre-defined parameters, and $\sigma(T_m)$ is the standard error of the energy function at the $m$th iteration. Within the first $n$ iterations, we use the statistical cooling method proposed by Aarts and Van Laarhoven (1985); while after the $n$th iteration, we adopt the geometric cooling strategy. Thirdly, notice that the neighborhood solution is simply generated by randomly selecting a node pair and changing its first front station, which means that the neighborhood solution makes just a little change to the current solution. We thus compute the energy value of a new solution by merely considering the changed portion rather than to re-compute it using Eq. (25) and (36) directly, resulting a significant reduction of solution time.

The modified model with penalty function is solved by the SA algorithm that has been implemented in C++ and run on the 3.30 GHz Intel Core i3 processor. We tested the algorithm with nine different instances (eight of them are from Section 4.3 and the remainder is constructed in this section) based on above algorithm rules. The results show that compared with Lingo, the RFFA can be solved by the SA with guaranteed solution quality and obviously improved time performance.

## 5. Case study

We conduct a case study based on the planned China rail network for the year 2020 and present the numerical results with visual illustrations from the proposed RFFA model and SA algorithm. We find that our method can effectively tackle the large-scale RFFA problem. The proposed SA algorithm is applied to the rail network for RFFA III.

To clearly show the tree-shaped path pattern of rail freight flow assignment on the real-world railway network, we present two typical links from the computational results for illustrative purposes. Fig. 5 shows the flow assignment pattern of all shipments destined for Shanghai region passing through the arc Suzhou→Shanghai. We observe that the freight flows that originated mostly from northwest and southwest China and converged to the tree root Shanghai form a tree-shaped path after merging several times.

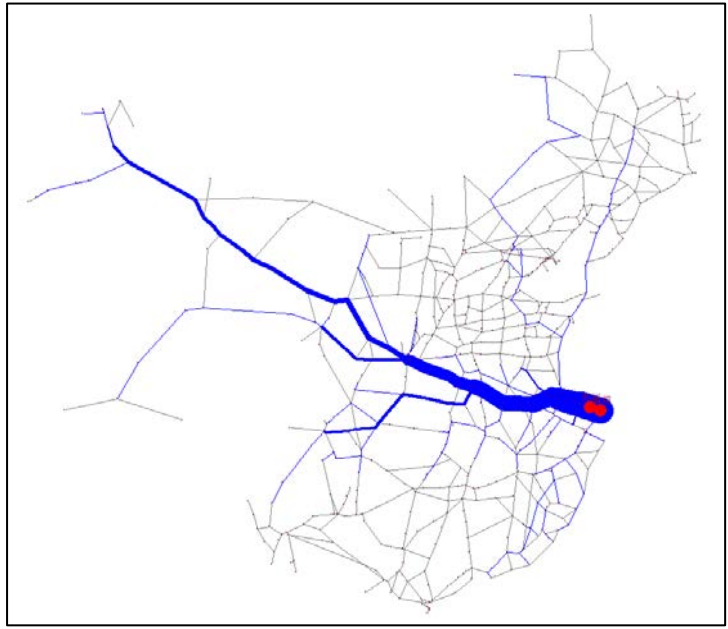

**Fig. 5.** Tree-shaped path of flows destined for Shanghai region.

To deepen our impression, without loss of generality, we show the flow assignment pattern of another typical arc of Qingyuan→Guangzhou, which is located in south China. Most of the O-D pairs passing through this arc come from the north, northwest, central and southwest of China. The path of these shipments also has an obvious tree-shaped characteristic that all flows converge at the tree root Guangzhou after merging several times (see Fig. 6).

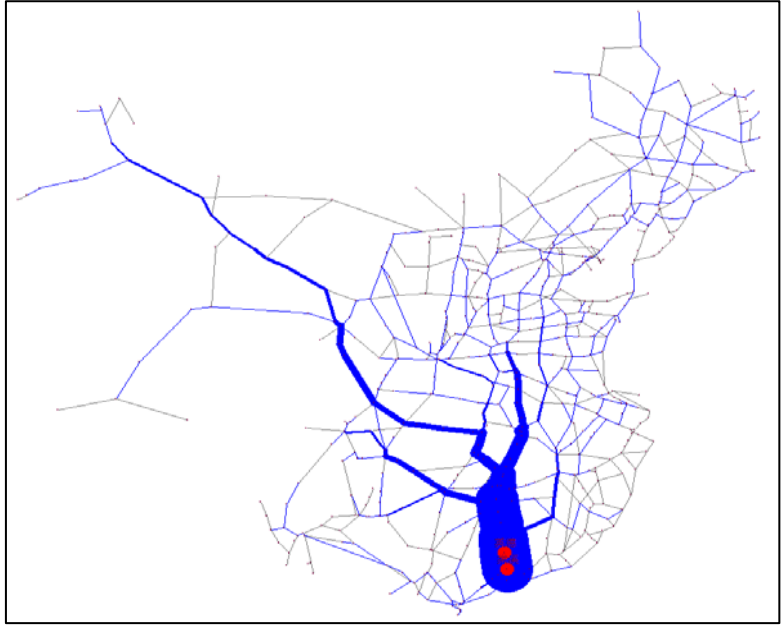

**Fig. 6.** Tree-shaped path of flows destined for Guangzhou region.

## 6. Conclusions

In this paper, we study the rail freight flow assignment problem with the tree-shaped path. Based on the tree feature, we develop a nonlinear programming model for the RFFA problem. The objective of the model is to minimize the total transport cost, and the constraints include the flow conservation constraint, tree-shaped path constraint and link capacity constraint. Because of the notably large number of quadratic terms generated from the tree-shaped constraints, the model is difficult to solve for a large-scale rail network. We thus introduce binary decision variables and modify the model into a nonlinear binary programming problem. In addition, we propose the penalty term method and the virtual arc method to address the potential problem of infeasible flows and analyze the complexity of the model.

The proposed model and solution approach have successfully solved the large-scale rail freight flow assignment flow problem with tree-shaped path constraints. To the best of our knowledge, previous studies have not addressed similar problems on such a large-scale rail network.

# References

Aarts, E.H.L., Van Laarhoven, P.J.M., 1985. Statistical cooling: a general approach to combinatorial optimization problems. *Philips Journal of Research* 40 (4), 193-226.

Černý, V., 1985. Thermodynamical approach to the traveling salesman problem: an efficient simulation algorithm. *Journal of Optimization Theory and Applications* 45 (1), 41-51.

Even, S., Itai, A., Shamir, A., 1975. On the complexity of time table and multi-commodity flow problems, *16th Annual Symposium on Foundations of Computer Science (SFCS 1975)*, USA, pp. 184-193.

Granville, V., Krivánek, M., Rasson, J.P., 1994. Simulated annealing: a proof of convergence. *IEEE Transactions on Pattern Analysis & Machine Intelligence* 16 (6), 652-656.

Kirkpatrick, S., Gelatt, C.D., Vecchi, M.P., 1983. Optimization by simulated annealing. *Science* 220 (4598), 671-680.

Kuby, M., Xu, Z., Xie, X., 2001. Railway network design with multiple project stages and time sequencing. *Journal of Geographical Systems* 3 (1), 25-47.

Li, P., 2000. Optimization of the train formation problem and extended models with variant parameters based on the genetic algorithm. Master thesis, Northern Jiaotong University, Beijing.

Lin, B.-L., Wang, Z.-M., Ji, L.-J., Tian, Y.-M., Zhou, G.-Q., 2012. Optimizing the freight train connection service network of a large-scale rail system. *Transportation Research Part B: Methodological* 46 (5), 649-667.

Lin, B., 2002. A tree-type car route guidance models for ITS, in: Wang, K.C.P., Xiao, G., Nie, L., Yang, H. (Eds.), *International Conference on Traffic and Transportation Studies (ICTTS) 2002*. American Society of Civil Engineers, Guilin, China, pp. 620-625.

Lin, B., Liu, C., Wang, H., Lin, R., 2017. Modeling the railway network design problem: A novel approach to considering carbon emissions reduction. *Transportation Research Part D: Transport and Environment* 56 (Supplement C), 95-109.

Metropolis, N., Rosenbluth, A.W., Rosenbluth, M.N., Teller, A.H., Teller, E., 1953. Equation of state calculations by fast computing machines. *The Journal of Chemical Physics* 21 (6), 1087-1092.

Selim, S.Z., Alsultan, K., 1991. A simulated annealing algorithm for the clustering problem. *Pattern Recognition* 24 (10), 1003-1008.

Tian, Y.M., Lin, B.L., Ji, L.J., 2011. Railway car flow distribution node-arc and arc-path models based on multi-commodity and virtual arc. *Journal of the China Railway Society* 33 (4), 7-12.

Wilson, R.J., 1996. *Introduction to graph theory (4th edition)*. Addison Wesley Longman Limited, Harlow, England.